\documentclass[11pt,a4paper]{article}
\usepackage[T1]{fontenc}
\usepackage{lmodern,microtype}
\usepackage[a4paper,left=28mm,right=28mm,top=24mm,bottom=27mm]{geometry}
\usepackage{amsmath,amssymb,amsthm}
\usepackage{tikz,flafter,needspace}
\usetikzlibrary{arrows.meta}
\usepackage[small,labelfont=bf,labelsep=period]{caption}
\usepackage[colorlinks=true,linkcolor=blue!45!black,citecolor=blue!45!black,urlcolor=blue!45!black]{hyperref}

\newtheorem{theorem}{Theorem}
\newtheorem{lemma}{Lemma}[section]
\newtheorem{corollary}{Corollary}[section]
\theoremstyle{definition}
\newtheorem{example}{Example}
\numberwithin{equation}{section}
\DeclareMathOperator{\conv}{conv}
\DeclareMathOperator{\suppv}{supp}
\definecolor{polygonblue}{HTML}{245369}
\definecolor{polygonfill}{HTML}{E1EFF2}
\definecolor{polygonrust}{HTML}{A04C32}
\definecolor{guidegrey}{HTML}{64737A}
\hypersetup{
  pdftitle={The eigenvalue region of doubly stochastic matrices near roots of unity},
  pdfauthor={Brecht Verbeken and Vincent Ginis}
}

\begin{document}
\newgeometry{left=28mm,right=28mm,top=18mm,bottom=23mm}
\begingroup
 \centering
 {\LARGE The eigenvalue region of doubly stochastic\\[2pt]
 matrices near roots of unity\par}
 \vspace{.5em}
 {\large Brecht Verbeken\textsuperscript{1,2} \qquad
 Vincent Ginis\textsuperscript{1,2,3}\par}
 \smallskip
 {\footnotesize
 \textsuperscript{1}Department of Business Technology and Operations, Data Analytics Laboratory,\\
 Vrije Universiteit Brussel (VUB), Brussels, Belgium\\
 \textsuperscript{2}imec-SMIT, Vrije Universiteit Brussel, Brussels, Belgium\\
 \textsuperscript{3}School of Engineering and Applied Sciences, Harvard University, Cambridge, MA, USA\par}
\endgroup

\begin{abstract}
Let \(D_N\) be the eigenvalue region of doubly stochastic matrices of order \(N\), and let
\(\Pi_k\) be the convex hull of the \(k\)th roots of unity. Near a root of exact order
\(q\le N\), \(D_N\) agrees with \(\Pi_M\), where \(M=q\lfloor N/q\rfloor\). Convex combinations
of powers of a cyclic shift realise \(\Pi_M\). For the side bounds we sum eigenvector
equations on the support of a limiting vector. The sum splits into cycles and paths.
For the two longest paths we use that the coordinates sum to zero. The union of the
polygons agrees with \(D_N\) in an annulus next to the unit circle. The same proof gives
the endpoint bound conjectured by Johnson, Newman and Spitkovsky.
\end{abstract}

\section{Introduction}\label{sec:intro}

A real \(N\times N\) matrix is doubly stochastic if its entries are nonnegative and every
row and column sums to one. For a permutation \(p\) of \(\{1,\ldots,N\}\), let \(P_p\)
have a one in column \(p(i)\) of each row \(i\) and zeros elsewhere. Thus
\((P_pu)_i=u_{p(i)}\) for column vectors \(u\in\mathbb C^N\).

In order five, let \(p\) interchange \(1\) and \(3\) and fix the other indices.
The matrix \(P_p\) sends \((1,0,-1,0,0)^\top\) to its negative and fixes the
four-dimensional subspace \(u_1=u_3\), so \(-1\) has algebraic multiplicity one.
For real matrices sufficiently close to \(P_p\), continuity of polynomial roots gives
one eigenvalue in \(|z+1|<1/2\), counted with multiplicity. It is real because a
nonreal eigenvalue would give a conjugate pair in the disk.

Let \(r\) interchange \(1\) with \(3\) and \(2\) with \(4\), and fix \(5\).
The eigenspace of \(P_r\) for \(-1\) consists of \((a,b,-a,-b,0)^\top\) with
\(a,b\in\mathbb C\). In Example~\ref{ex:square} we combine \(P_r\) with the matrices
for \((1\,2\,3\,4)\) and its inverse. This gives the eigenvalues \(-1+t(1\pm i)\),
\(0\le t\le1\), along two sides of a square at \(-1\). Do those sides bound the
eigenvalues near \(-1\) for every doubly stochastic matrix of order five?

Let \(\sigma(A)\) denote the set of eigenvalues of \(A\). Define
\[
 D_N=\bigcup_{\substack{A\ \mathrm{doubly\ stochastic}\\ A\in\mathbb R^{N\times N}}}
              \sigma(A),
 \qquad
 \Pi_k=\conv\{e^{2\pi i j/k}:0\le j<k\},
 \qquad
 U_N=\bigcup_{k=1}^N\Pi_k .
\]
Here \(\Pi_1\) is the point \(\{1\}\) and \(\Pi_2\) is the segment \([-1,1]\).
To realise a point of \(\Pi_k\), take a convex combination of powers of a cyclic
shift on \(k\) coordinates and add \(N-k\) fixed coordinates. In particular,
\(U_N\subseteq D_N\). Perfect and Mirsky proved equality through order three
and asked whether the inclusion could always be reversed \cite{pm}. Levick,
Pereira and Kribs proved equality in order four \cite{lpk}, and Mashreghi and
Rivard found an eigenvalue outside \(U_5\) \cite{mr}.

If \(q\le N\) is the least positive integer with \(\omega^q=1\), the polygons \(\Pi_k\) that
contain \(\omega\) are exactly those with \(q\mid k\). Does the polygon with the largest such
\(k\) already give all of \(D_N\) near \(\omega\)?

\begin{theorem}\label{thm:local}
Let \(\omega\) be a root of unity of exact order \(q\le N\), so \(q\) is the least positive integer for
which \(\omega^q=1\). Let \(M=q\lfloor N/q\rfloor\). There exists \(\varepsilon>0\) such that
\[
 D_N\cap B(\omega,\varepsilon)
   = U_N\cap B(\omega,\varepsilon)
   = \Pi_M\cap B(\omega,\varepsilon),
\]
where \(B(\omega,\varepsilon)=\{z\in\mathbb C:|z-\omega|<\varepsilon\}\).
For each fixed \(N\) there is also an \(r_N\in[0,1)\) such that
\[
 D_N\cap\{|z|>r_N\}=U_N\cap\{|z|>r_N\}.
\]
\end{theorem}
\restoregeometry

Convex combinations of the powers of a single permutation matrix have their
eigenvalues in \(U_N\). Varying that matrix over all permutation matrices of order \(N\)
gives all of \(U_N\). This follows from the description for finite abelian matrix
groups by Jankowski, Johnson and Lim \cite[Corollary~3.5]{jjl}. For the same class
of convex combinations, Casta\~neda Betancourt, Morneau-Gu\'erin and Singh
determine the least order that realises a given eigenvalue
\cite[Section~4.2]{cms}.

Harlev, Johnson and Lim conjectured that every boundary point of \(D_N\) is an
eigenvalue of a matrix on a segment between two permutation matrices
\cite[Conjecture~2.7]{hjl}. At a fixed permutation matrix \(P\), Johnson, Newman
and Spitkovsky studied a branch issuing from an eigenvalue \(\omega\) along
\(A(t)=(1-t)P+tR\), with \(R\) doubly stochastic. For \(\omega\) of exact order
\(q\), the bound they propose uses the total length of all cycles of the permutation for \(P\)
whose lengths are divisible by \(q\). We write
\[
 L=\sum_{\substack{C\text{ a cycle associated with }P\\ q\mid |C|}} |C|.
\]
The conjecture is that a continuous eigenvalue branch with \(z(0)=\omega\)
lies in \(\Pi_L\) for all sufficiently small \(t\ge0\)
\cite[Section~3, just before Theorem~3.2]{jns}. Their Theorem~3.2 treats the
case of a single cycle whose length is divisible by \(q\).
Corollary~\ref{cor:fixed} gives this bound for all eigenvalues near \(\omega\)
of doubly stochastic matrices sufficiently close to \(P\).

The count \(L\) bounds the number of nonzero coordinates of an eigenvector of
\(P\) for \(\omega\). Indeed, suppose \(A_j u_j=z_j u_j\), \(A_j\to P\) and
\(z_j\to\omega\). A subsequential limit \(v\) of unit eigenvectors satisfies
\(Pv=\omega v\). Going once around a cycle of length \(\ell\) gives \(v_i=\omega^\ell v_i\) at
each index \(i\) on that cycle. If \(v_i\ne0\), then \(q\mid\ell\).
Thus the support \(S=\{i:v_i\ne0\}\) satisfies \(|S|\le L\le M\), since
\(L\) is a multiple of \(q\) at most \(N\). The opening
matrices \(P_p\) and \(P_r\) give \(L=2\) and \(L=4\) at \(-1\) in order five,
where \(M=4\).

For the full region we only require \(z_j\to\omega\). Take a subsequence
on which \(A_j\to A\) and the unit eigenvectors \(u_j\to v\). Then
\(Av=\omega v\) and \(\|v\|_2=1\). In a decomposition of \(A\) into
permutation matrices, every term with positive
coefficient sends \(v\) to \(\omega v\), by equality in the norm bound.
Following the cycles of one such permutation shows that \(|S|\) is a multiple
of \(q\). Since \(|S|\le N\), the size bound is \(|S|\le M\).
For \(q\ge2\), Lemma~\ref{lem:stable} gives the two sides of \(\Pi_M\) using
this bound if \(M\ge3\). For matrices approaching \(P\), we have the bound
\(|S|\le L\), so the same lemma gives the sides of \(\Pi_L\) if \(L\ge3\).
For \(L=q=2\), the eigenvalue is simple as in the opening example. For
\(q=1\), both \(L\) and \(M\) equal \(N\).

For an eigenpair \(Av=zv\), each row equation expresses \(zv_i\) as a convex
combination of the coordinates of \(v\). Hence \(zK\subseteq K\) for
\(K=\conv\{v_1,\ldots,v_N\}\). Dmitriev and Dynkin used this invariant polygon
for row-stochastic matrices \cite{dd}. Karpelevi\v{c} determined their eigenvalue
region, whose boundary joins roots of unity by arcs \cite{karp,vgkarp}.
For doubly stochastic matrices, Levick, Pereira and Kribs used the column sums
to obtain the same invariance for the hull of the sums \(v_i+v_j\), \(i<j\)
\cite[Corollary~2]{lpk}. Bouthat, Mashreghi and Morneau-Gu\'erin consider hulls of
sums over any number of distinct indices. They announce \(z\in\Pi_N\) when
each coordinate lies outside the convex hull of the others
\cite[slides~28--31]{bmm}. In our order-five example, the eigenvector
\((1,1,-1,-1,0)^\top\) for \(-1\) has equal coordinates at distinct indices.

For a nearby eigenpair \(Au=zu\), multiply the eigenvector equation in row
\(i\in S\) by \(\bar u_i\) and add over \(S\). On the right we get an average
over choices of distinct columns \(p(i)\) for these rows, one choice for
each permutation \(p\). At the limit, compare each choice with a permutation
sending \(v\) to \(\omega v\). The sum then separates into cycles and paths
on \(S\). We estimate cycles by their roots of unity and paths by using the
zero sum of the limiting coordinates in the longest cases. At equality the
permutation preserves \(S\), and its spectrum there gives the side bound for
nearby vectors.

\section{Cyclic shifts and the two sides of a polygon}\label{sec:polygons}

\subsection{Realising the polygons}

Let \(C_k\) be the cyclic shift \((C_kw)_j=w_{j+1}\) with \(w_{k+1}=w_1\).
To realise a convex combination of \(k\)th roots of unity, we use one
eigenvector for all the powers of \(C_k\). For \(\zeta=e^{2\pi i/k}\) the vector
\(w=(1,\zeta,\ldots,\zeta^{k-1})^\top\) satisfies \(C_kw=\zeta w\), since
\(\zeta^k=1\) closes the shift at the last coordinate. If
\(a_0,\ldots,a_{k-1}\) are nonnegative and sum to one,
\[
 \left(\sum_{j=0}^{k-1}a_j C_k^j\right)w
 =\left(\sum_{j=0}^{k-1}a_j\zeta^j\right)w .
\]
Varying the coefficients gives every eigenvalue in \(\Pi_k\), with a
doubly stochastic matrix on the left. Adding an identity block of order
\(N-k\) when \(k<N\) leaves the eigenvalue unchanged on the vector obtained
by appending zeros to \(w\). This proves \(U_N\subseteq D_N\).

\begin{example}[The sides at \(-1\) in order five]\label{ex:square}
Let \(v=(1,i,-1,-i,0)^\top\), and write \(\oplus\) for a block diagonal
sum. Starting from \(C_4^2\oplus[1]\), move towards \(C_4\oplus[1]\) by taking
the following convex combination for \(0\le t\le1\).
\[
 A(t)=\big((1-t)C_4^2+tC_4\big)\oplus[1]
     =\begin{pmatrix}
        0&t&1-t&0\\
        0&0&t&1-t\\
        1-t&0&0&t\\
        t&1-t&0&0
       \end{pmatrix}\oplus[1].
\]
On the first four coordinates of \(v\), the matrices \(C_4\) and \(C_4^2\)
act by multiplication by \(i\) and \(i^2=-1\), respectively. Thus
\(A(t)v=(-1+t(1+i))v\), and the eigenvalue moves from \(-1\) to \(i\).
Replacing \(C_4\) by \(C_4^3\) makes it move from \(-1\) to \(-i\).
Both matrix families start at \(P_r=C_4^2\oplus[1]\).
\end{example}

\subsection{The sides at a root of unity}

Let \(\omega\) have exact order \(q\le N\). Project each \(k\)th root onto the
radius through \(\omega\) by taking \(\Re(\bar\omega z)\). The value is at
most \(1\), and it equals \(1\) only at \(\omega\). A convex combination
can reach \(\omega\) only if \(\omega\) itself is a \(k\)th root. Thus
\(\omega\in\Pi_k\) exactly when \(q\mid k\). Each other polygon is compact
and is at positive distance from \(\omega\).

Take \(k\ge3\) with \(q\mid k\) and put \(\alpha=\pi/k\). The sides at
\(\omega\) join it to \(\omega e^{2i\alpha}\) and \(\omega e^{-2i\alpha}\).
The outward unit normals point halfway between the respective radii, so
they are \(\eta_\pm=\omega e^{\pm i\alpha}\). The projection of \(z\)
onto a normal is \(\Re(\bar\eta_\pm z)\). Both endpoints of the corresponding
side have projection \(\cos\alpha\), giving
\begin{equation}\label{eq:side}
 \Re(\bar\eta_\pm z)\le \cos\alpha .
\end{equation}
The other sides give strict inequalities at \(\omega\). They therefore hold
throughout a small disk about \(\omega\), where only \eqref{eq:side} is needed.

It is convenient to measure an inward distance \(x\) and a signed distance \(y\)
perpendicular to the radius, writing \(z=\omega(1-x+iy)\). The two projections are
\((1-x)\cos\alpha\pm y\sin\alpha\). Thus both inequalities together read
\begin{equation}\label{eq:wedge}
 x\ge |y|\tan(\pi/k).
\end{equation}
For a given \(x\), the range of \(y\) grows with \(k\). So near \(\omega\) the largest
multiple \(M\) of \(q\) with \(M\le N\) gives all of \(U_N\). At an endpoint of
\(\Pi_2\) the local set is the inward radial segment, and \(\Pi_1\) contributes just
the point \(1\). Both are included in the same local description.

\begin{figure}[!htbp]
 \centering
 \begin{tikzpicture}[x=1.72cm,y=1.72cm,>=Latex,font=\small]
  \fill[polygonfill] (0:1)--(90:1)--(180:1)--(270:1)--cycle;
  \draw[guidegrey!65,thin] (0,0) circle[radius=1];
  \draw[guidegrey!55,thin] (-1.14,0)--(1.14,0) (0,-1.12)--(0,1.12);
  \draw[polygonblue,thick] (0:1)--(90:1)--(180:1)--(270:1)--cycle;
  \draw[polygonrust,semithick,dashed]
     (0:1)--(72:1)--(144:1)--(216:1)--(288:1)--cycle;
  \fill[polygonblue] (180:1) circle[radius=1.45pt] node[left=2pt] {\(-1\)};
  \node[polygonblue,above left] at (0,1) {\(i\)};
  \node[polygonblue,below left] at (0,-1) {\(-i\)};
  \node[polygonblue] at (0.28,0.34) {\(\Pi_4\)};
  \node[polygonrust,right] at (0.74,0.53) {\(\Pi_5\)};
  \begin{scope}[xshift=6.55cm,x=3.30cm,y=3.30cm]
   \fill[polygonfill] (0,0)--(.46,.46)--(.46,-.46)--cycle;
   \draw[guidegrey!55,thin] (-.14,0)--(.53,0);
   \draw[polygonblue,thick] (.46,.46)--(0,0)--(.46,-.46);
   \draw[polygonrust,semithick,dashed] (.191,-.49)--(.191,.49);
   \draw[polygonblue,semithick,->] (.115,.115)--(.233,.233);
   \draw[polygonblue,semithick,->] (.115,-.115)--(.233,-.233);
   \fill[polygonblue] (0,0) circle[radius=1.45pt] node[left=2pt] {\(-1\)};
   \node[polygonrust,above] at (.191,.49) {\(\Pi_5\)};
   \node[polygonblue,right] at (.46,.27) {\(\Pi_4\)};
  \end{scope}
 \end{tikzpicture}
 \caption{The square \(\Pi_4\) and pentagon \(\Pi_5\) in order five.
 The right panel shows the two sides at \(-1\) traced by the families
 in Example~\ref{ex:square}. The dashed edge of \(\Pi_5\) is at positive distance from \(-1\).}
 \label{fig:square}
\end{figure}
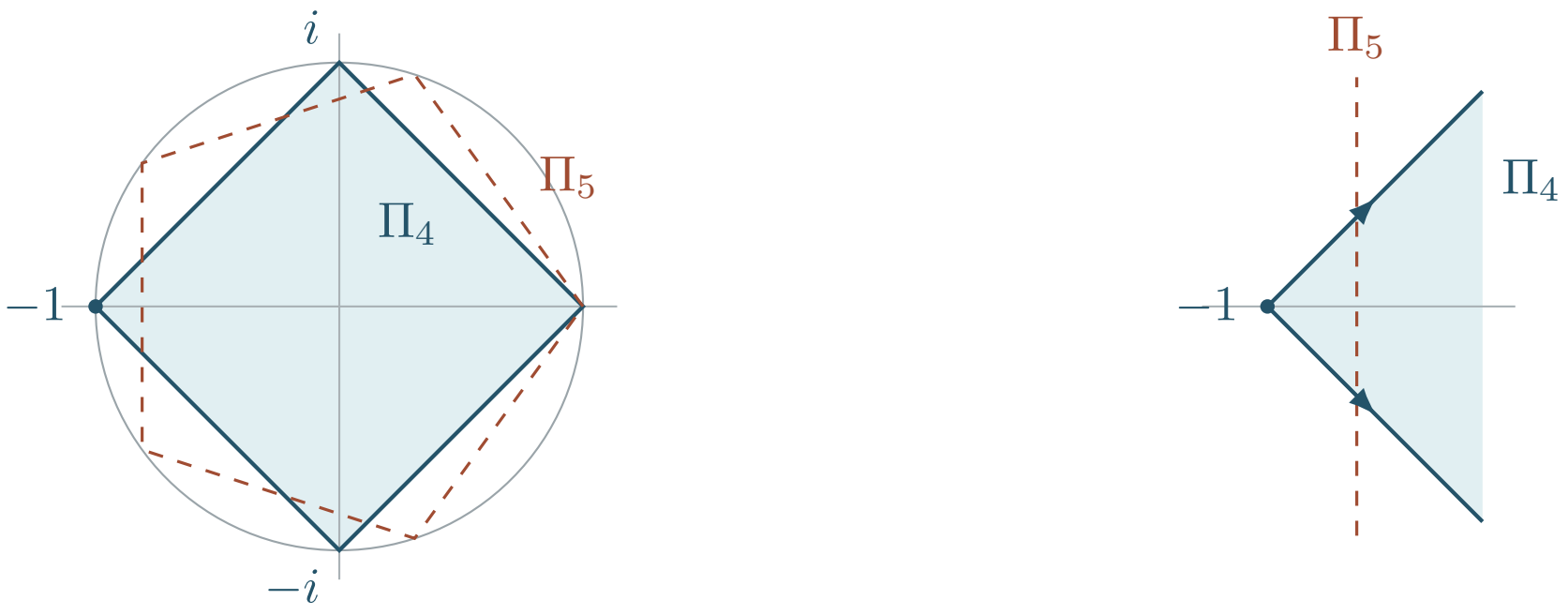

For \(N=8\), the root \(\omega=e^{2\pi i/3}\) is a vertex of the triangle
\(\Pi_3\) and the hexagon \(\Pi_6\). At inward distance \(x\), they allow
\(|y|\le x/\sqrt3\) and \(|y|\le\sqrt3\,x\), respectively. The hexagon
has the wider angle. The octagon is at positive distance from
\(\omega\), so the side inequalities to prove for \(D_8\) are those of
\(\Pi_6\).

\section{The eigenvector equation on selected rows}\label{sec:assignment}

A permutation matrix \(P_p\) selects one column \(p(i)\) in each row \(i\).
The map \(p\) is one-to-one, so two rows cannot select the same column.
We use the following decomposition to retain this restriction when we
work with a doubly stochastic matrix.

\begin{lemma}[Birkhoff--von Neumann]\label{lem:birkhoff}
Every doubly stochastic matrix is a convex combination of permutation matrices of the same order.
\end{lemma}

\begin{proof}
For a doubly stochastic matrix \(A\), join each row to the columns in which its entries
are positive. If \(I\) is a set of rows and \(\Gamma(I)\) is the set of their neighbouring
columns, the entries in the rows of \(I\) have sum \(|I|\). Each column of \(\Gamma(I)\)
contributes at most one, so \(|I|\le|\Gamma(I)|\).

We need one positive entry per row with distinct columns. Choose a matching
of maximum size, where a matching is a set of edges meeting each row and
column at most once. Suppose a row is unmatched. From it follow unmatched
edges to columns and matched edges back to rows. If we reach an unmatched
column, a shortest path to it starts and ends with an unmatched edge.
Swapping the matched and unmatched edges along this path adds one edge to
the matching. Thus every reached column is matched to a reached row.
The starting row has no partner, so there are fewer reached columns than
reached rows. Every neighbour of a reached row is reached. An unmatched
edge can be followed from the row, and a matched edge is how the row was
reached. This contradicts \(|I|\le|\Gamma(I)|\) for the reached rows \(I\).

The chosen entries now give a permutation matrix \(P\). Let \(a\) be the least of
these positive entries. If \(a=1\), then \(A=P\). If \(a<1\), the matrix
\((A-aP)/(1-a)\) is again doubly stochastic and has fewer positive entries. Apply the
same argument to it. Each subtraction removes a positive entry, so the process
terminates and expresses \(A\) as a convex combination of the chosen matrices.
\end{proof}

\Needspace{7\baselineskip}
The numerical range of a square matrix \(B\) is the set of values
\(u^*Bu\) for \(\|u\|_2=1\). Here \(u^*\) is the conjugate transpose and
\(\|\cdot\|_2\) is the Euclidean norm. Unit eigenvectors give the
eigenvalues among these values. If \(B\) is normal, write
\(u=\sum_j a_jw_j\) in an orthonormal eigenbasis with \(Bw_j=\lambda_j w_j\).
Then \(u^*Bu=\sum_j |a_j|^2\lambda_j\). The weights \(|a_j|^2\) range over
all nonnegative numbers of sum one, so the numerical range of \(B\) is
the convex hull of its eigenvalues.

Permutation matrices are unitary. A cycle of length \(\ell\) gives
a block \(C_\ell\), which has the \(\ell\)th roots of unity as its
eigenvalues by the calculation in Section~\ref{sec:polygons}. Thus
the numerical range for each permutation matrix lies in \(\conv U_N\).
For a unit vector \(u\) and a doubly stochastic \(A\), Lemma~\ref{lem:birkhoff}
expresses \(u^*Au\) as a convex combination of those values. This
gives \(D_N\subseteq\conv U_N\).

Consider the side from \(-1\) to \(i\) in Figure~\ref{fig:square}.
Its line has projection \(\cos(\pi/4)\) onto the normal \(e^{3\pi i/4}\).
The fifth root \(e^{4\pi i/5}\) projects to \(\cos(\pi/20)\), which is
larger. The segment in \(\conv U_5\) joining \(-1\) to this root
therefore contains points beyond the side arbitrarily close to
\(-1\). To obtain the local side bound we sum over selected rows.

Let \(S\subseteq\{1,\ldots,N\}\) be the rows we keep, with total weight
\(W_S(u)=\sum_{i\in S}|u_i|^2\) for \(u\in\mathbb C^N\). For
\(|\eta|=1\), we write the projection of the sum for \(p\) as
\begin{equation}\label{eq:score}
 F_\eta(p,u)=\Re\left(\bar\eta\sum_{i\in S}\bar u_i u_{p(i)}\right).
\end{equation}
Only the assignment of columns to the retained rows enters the sum.
These columns are distinct. Conversely, one can complete any injection
from \(S\) into \(\{1,\ldots,N\}\) to a permutation by assigning the
unused columns to the other rows. Equal coordinates can still have
different column indices.

Take an eigenpair \(Au=zu\). Multiplying row \(i\) by \(\bar u_i\) and
summing over \(S\) gives \(zW_S(u)=\sum_{i\in S}\bar u_i(Au)_i\).
By Lemma~\ref{lem:birkhoff}, the projection on the right is an
average of \(F_\eta(p,u)\) with the coefficients of the permutation
decomposition. It is bounded by the largest term, so
\begin{equation}\label{eq:scorebound}
 W_S(u)\Re(\bar\eta z)\le\max_p F_\eta(p,u).
\end{equation}
Here only \(p\) varies in the maximum. When \(W_S(u)>0\), a bound on
the weighted sums gives a bound on the projection of \(z\).

\begin{example}[Two rows competing for a column]\label{ex:distinct}
Take \(u=(1,1,-2,0)^\top\), \(S=\{1,2\}\) and \(\eta=-1\). The projection
in \eqref{eq:score} is \(-u_{p(1)}-u_{p(2)}\). Both rows would gain from
column \(3\), and giving it to both would give \(4\). With distinct
columns the best choice is columns \(3\) and \(4\), giving \(2\).
The unprojected sum takes the values \(-2,-1,1,2\), whose convex hull
is \([-2,2]\). Since \(W_S(u)=2\), the bound in direction \(\eta=-1\)
after dividing by the weight is \(1\).
\end{example}

\Needspace{11\baselineskip}
\section{The support of an eigenvector at a root of unity}\label{sec:support}

By Lemma~\ref{lem:birkhoff}, \(Au\) is a convex combination of vectors
of norm \(\|u\|_2\). Thus \(\|Au\|_2\le\|u\|_2\) and \(D_N\) lies in
the closed unit disk. At a point on the circle the norm bound is an
equality. This will tell us which rows to keep in \eqref{eq:scorebound}.

\begin{lemma}\label{lem:support}
Suppose \(Av=\omega v\), with \(A\) doubly stochastic, \(v\ne0\) and
\(|\omega|=1\). In any decomposition of \(A\) as a convex combination
of permutation matrices, each matrix with positive coefficient sends
\(v\) to \(\omega v\). In particular, \(\omega\) is a root of unity of
exact order \(q\le N\).

Write \(S=\suppv v=\{i:v_i\ne0\}\) and \(m=|S|\). If \(q\ge2\), then
\(q\mid m\), the sum \(\sum_{i\in S}v_i\) is zero, and at least \(q\)
coordinates of \(v\) have largest modulus.
\end{lemma}

\begin{proof}
Take a decomposition \(A=\sum_p a_p P_p\) from Lemma~\ref{lem:birkhoff}.
To see what equality of norms implies, compare each \(P_pv\) with
\(\omega v\). Since \(\|P_pv\|_2=\|v\|_2\) and
\(v^*Av=\omega\|v\|_2^2\),
\[
 \sum_p a_p\|P_pv-\omega v\|_2^2
   =2\|v\|_2^2-2\Re(\bar\omega\,v^*Av)=0 .
\]
Each squared distance with \(a_p>0\) vanishes. Choose one permutation
\(r\) with positive coefficient. The equations \(v_{r(i)}=\omega v_i\)
show that a cycle lies in \(S\) as soon as it contains a nonzero coordinate.
On such a cycle of length \(\ell\), after \(\ell\) steps we have
\(v_i=\omega^\ell v_i\), hence \(\omega^\ell=1\). The exact order \(q\)
exists, and it divides the length of every cycle contained in \(S\).
In particular \(q\le N\) and \(q\mid |S|\).

When \(q\ge2\), divide each of these cycles into consecutive groups
of \(q\) indices. The coordinates in a group have the form
\(v_i,\omega v_i,\ldots,\omega^{q-1}v_i\). They sum to zero, and all
\(q\) have the same modulus. Summing over the groups gives zero,
and the group containing a coordinate of largest modulus has \(q\)
coordinates of that modulus.
\end{proof}

For \(C_4^2\oplus[1]\) from Example~\ref{ex:square}, the vectors
\((1,i,-1,-i,0)^\top\), \((1,1,-1,-1,0)^\top\) and
\((1,0,-1,0,0)^\top\) all belong to the \(-1\)-eigenspace. The first
two use all four indices of the two cycles, even though the coordinates
repeat in the second vector. The third uses just the cycle on
\(\{1,3\}\).

Let \(A_j u^{(j)}=z_j u^{(j)}\), where \(z_j\to\omega\) and
\(\|u^{(j)}\|_2=1\). The set of doubly stochastic matrices and the unit sphere are
compact, so after taking a subsequence we have \(A_j\to A\) and \(u^{(j)}\to v\).
Then \(Av=\omega v\) and \(\|v\|_2=1\). In \eqref{eq:scorebound} use
\(S=\suppv v\) for all \(u^{(j)}\), even when \(u^{(j)}\) has
nonzero coordinates outside \(S\). We have \(W_S(v)=1\) and
\(W_S(u^{(j)})\to1\). The next example explains why we keep this
set of rows.

\begin{example}[The numerical range near a limiting eigenvector]\label{ex:small-coordinates}
In order seven take \(P=C_2\oplus C_5\), \(\omega=-1\) and \(\lambda=e^{6\pi i/5}\).
Let \(v=(1,-1,0,0,0,0,0)^\top/\sqrt2\) and
\(w=(0,0,1,\lambda,\ldots,\lambda^4)^\top/\sqrt5\), so \(Pv=-v\) and
\(Pw=\lambda w\). The side of \(\Pi_6\) from \(-1\) to \(-e^{i\pi/3}\) has
normal \(\eta=-e^{i\pi/6}\). The projection of \(\lambda\) onto it is
\(\cos(\pi/30)>\cos(\pi/6)\).

In \eqref{eq:scorebound} the eigenvector of \(A\) is tested against
individual permutation matrices. It may have components in several
eigenspaces of the matrix against which it is tested. As a test vector for \(P\),
take \(u=v+sw\) with \(s>0\). The disjoint supports of the two unit
vectors give
\[
 \frac{u^*Pu}{\|u\|_2^2}=\frac{-1+s^2\lambda}{1+s^2}.
\]
The \(s^2\lambda\) term sends the quotient beyond the side for every
\(s>0\), however small. On the support \(S=\{1,2\}\) of \(v\),
both \(u\) and \(Pu\) agree with \(v\) and \(Pv\), respectively.
Thus the sum over \(S\) stays at \(-1\) and \(W_S(u)=1\).
\end{example}

\section{Cycles and paths on the support}\label{sec:paths}

Let \(P_rv=\omega v\), where \(v\ne0\) and \(\omega\) has exact order
\(q\ge2\). Write \(S=\suppv v\) and \(m=|S|\). To estimate the side
contributions we only need an upper bound \(M\) on \(m\) with
\(q\mid M\) and \(M\ge3\). For Theorem~\ref{thm:local} we take the
size bound \(M=q\lfloor N/q\rfloor\). For Corollary~\ref{cor:fixed} we
take \(M=L\) when \(L\ge3\), using the cycles of the fixed permutation.
For the chosen \(M\), put
\(\alpha=\pi/M\) and \(\xi=e^{2i\alpha}\). We start with the side
from \(\omega\) to \(\omega\xi\), of normal \(\eta=\omega e^{i\alpha}\).

For a permutation \(p\) in \eqref{eq:scorebound}, compare the column
\(p(i)\) with the columns of the fixed permutation \(r\).
There is a unique index \(j\) with \(r(j)=p(i)\). This defines the
permutation \(T=r^{-1}\circ p\) by \(T(i)=j\). Now
\(v_{p(i)}=\omega v_{T(i)}\), and the factor \(\omega\) cancels
with the factor \(\bar\omega\) in \(\bar\eta\). Hence
\begin{equation}\label{eq:relative}
 F_\eta(p,v)=\Re\left(e^{-i\alpha}\sum_{i\in S}\bar v_i v_{T(i)}\right).
\end{equation}
In \eqref{eq:relative} a term survives precisely when both \(i\) and
\(T(i)\) lie in \(S\). Keep the arrows \(i\to T(i)\) for these terms.
A cycle wholly inside \(S\) remains a cycle, including a loop at a
fixed point. Removing the indices outside \(S\) cuts the other
cycles into directed paths. A vertex incident to no retained arrow
counts as a path on one vertex.

\begin{figure}[!htbp]
 \centering
 \begin{tikzpicture}[
  x=1.05cm,y=1.08cm,>=Latex,
  vertex/.style={circle,draw=polygonblue,fill=polygonfill,thick,minimum size=6.7mm,inner sep=0pt},
  outside/.style={vertex,draw=guidegrey,fill=white,dashed},
  arrow/.style={->,draw=polygonblue,semithick},
  omitted/.style={->,draw=guidegrey,semithick,dashed},
  note/.style={font=\small}
 ]
  \node[vertex] (n1) at (0,0) {\(1\)};
  \node[vertex] (n2) at (1,1) {\(2\)};
  \node[vertex] (n3) at (2,0) {\(3\)};
  \draw[arrow] (n1)--(n2);
  \draw[arrow] (n2)--(n3);
  \draw[arrow] (n3)--(n1);
  \node[note,polygonblue,below] at (1,-.32) {Cycle on \(S\)};

  \node[vertex] (n4) at (3.5,.42) {\(4\)};
  \node[vertex] (n5) at (4.85,.42) {\(5\)};
  \node[vertex] (n6) at (6.2,.42) {\(6\)};
  \node[outside] (n7) at (7.55,.42) {\(7\)};
  \draw[arrow] (n4)--(n5);
  \draw[arrow] (n5)--(n6);
  \draw[omitted] (n6)--(n7);
  \draw[omitted] (n7) to[bend left=37] (n4);
  \node[note,polygonblue,above] at (4.85,.71) {Path on \(S\)};
  \node[note,guidegrey,above] at (7.55,.71) {Outside \(S\)};
 \end{tikzpicture}
 \caption{For \(T=(1\,2\,3)(4\,5\,6\,7)(8)\) and \(S=\{1,\ldots,6\}\),
 only the solid arrows contribute to \eqref{eq:relative}. Removing \(7\)
 leaves the path \(4\to5\to6\), and the cycle on \(\{1,2,3\}\) stays intact.}
 \label{fig:components}
\end{figure}
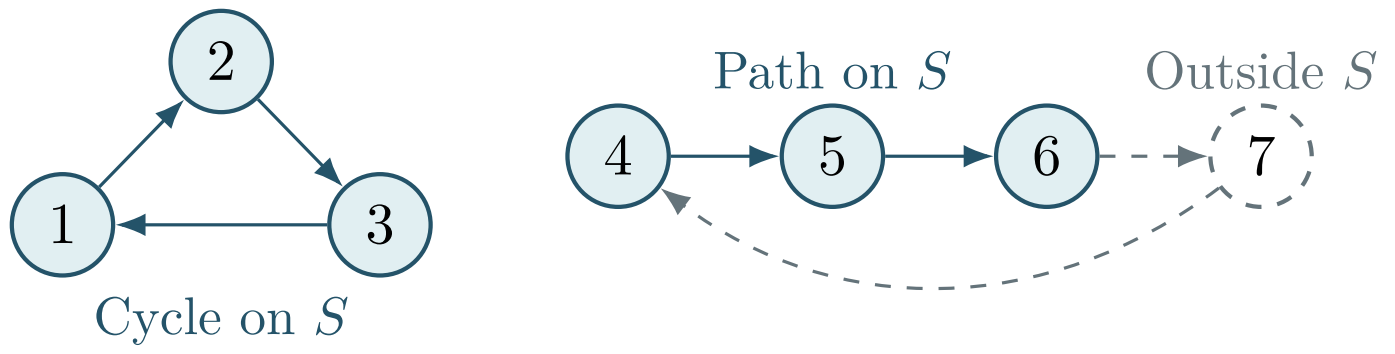

Take \(v=(1,-1,1,-1,1,-1,0,0)^\top\) and \(r=(1\,2)(3\,4)(5\,6)\),
fixing \(7\) and \(8\), so \(P_rv=-v\). For \(T\) from
Figure~\ref{fig:components}, let \(p=r\circ T\). The terms under the
projection in \eqref{eq:relative} sum to
\[
 (\bar v_1v_2+\bar v_2v_3+\bar v_3v_1)
  +(\bar v_4v_5+\bar v_5v_6)=-1-2 .
\]
The first sum follows the cycle and the second follows the
path. On each component we bound the projection in direction
\(e^{i\alpha}\) by \(\cos\alpha\) times the sum of squared moduli
on that component. Since each index in \(S\) occurs in one
component, adding the bounds gives
\(F_\eta(p,v)\le\cos\alpha\,W_S(v)\).

\subsection{The contribution of a cycle}\label{sec:cycle-comparison}

Write the coordinates around a cycle of length \(\ell\) in their
cyclic order in \(w\). The arrows give the sum \(w^*C_\ell w\), and
the numerical range of \(C_\ell\) is \(\Pi_\ell\). The side from
\(1\) to \(\xi\) cuts off the open arc of angle \(2\pi/M\) between
its endpoints. Since \(\ell\le M\), no \(\ell\)th root lies in
this arc. Thus
\begin{equation}\label{eq:cycle-comparison}
 \Re(e^{-i\alpha}w^*C_\ell w)\le \cos\alpha\,\|w\|_2^2 .
\end{equation}

In the cyclic shift eigenbasis the projection is a sum of the
projections of roots with nonnegative weights. Equality in
\eqref{eq:cycle-comparison} allows nonzero weights only at
roots on the side line. For \(\ell<M\), only \(1\) lies on the
line, so \(w\) is constant. For \(\ell=M\), the roots on the line
are \(1\) and \(\xi\), so \(w\) is a linear combination of the
constant vector and \((1,\xi,\ldots,\xi^{M-1})^\top\).
A cycle of length \(M\) in \eqref{eq:relative} uses every index of
\(S\). Its coordinates sum to zero by Lemma~\ref{lem:support},
so the constant term vanishes. The equality vector is then a
nonzero multiple of \((1,\xi,\ldots,\xi^{M-1})^\top\), whose
coordinates are all distinct.

\Needspace{8\baselineskip}
\subsection{Closing a path with a zero coordinate}

For a path of \(\ell\) vertices, write \(w=(w_1,\ldots,w_\ell)^\top\)
in the order of the arrows. To use the cycle bound, insert a zero
between \(w_\ell\) and \(w_1\). In the vector
\(\tilde w=(0,w_1,\ldots,w_\ell)^\top\) the two new arrows contribute
\(\overline{0}w_1\) and \(\bar w_\ell 0\). Hence
\begin{equation}\label{eq:path-completion}
 \tilde w^*C_{\ell+1}\tilde w=\sum_{j=1}^{\ell-1}\bar w_jw_{j+1},
 \qquad
 \|\tilde w\|_2^2=\|w\|_2^2 .
\end{equation}
If \(\ell\le M-2\), this cycle has length less than \(M\).
Equality in \eqref{eq:cycle-comparison} would make \(\tilde w\)
constant, which is impossible because it has a zero coordinate
and \(w\ne0\). Thus the path inequality is strict.

\subsection{A path on \texorpdfstring{\(M-1\)}{M-1} vertices}

For \(\ell=M-1\), the added zero gives a cycle of length \(M\).
The cycle estimate still applies, with two eigenspaces available
at equality. Number the coordinates from \(0\) to \(M-1\), starting
with the zero. At equality \(\tilde w_j=b+a\xi^j\), and
\(\tilde w_0=0\) gives \(b=-a\). Thus
\[
 w_j=a(\xi^j-1)\quad(1\le j\le M-1),\qquad a\ne0 .
\]
All these coordinates have modulus at most \(2|a|\) and
\(\sum_{j=0}^{M-1}\xi^j=0\) gives their sum \(-Ma\).
Since \(M-1\le m\le M\) and \(m,M\) are multiples of \(q\ge2\),
we have \(m=M\). Only one coordinate of \(S\) is left outside the
path. The zero sum on \(S\) forces that coordinate to be \(Ma\).
Its modulus exceeds \(2|a|\) as \(M>2\). This contradicts
Lemma~\ref{lem:support}, which gives at least \(q\) indices with
largest modulus. Thus the inequality is strict.

\begin{figure}[!htbp]
 \centering
 \begin{tikzpicture}[x=1.20cm,y=1.20cm,>=Latex,font=\small]
  \draw[guidegrey!60,->] (-2.55,0)--(4.75,0) node[right] {\(\Re\)};
  \draw[guidegrey!60,->] (0,-1.38)--(0,1.36) node[above] {\(\Im\)};
  \draw[polygonblue,semithick,->] (-1,1)--(-2,0);
  \draw[polygonblue,semithick,->] (-2,0)--(-1,-1);
  \fill[polygonblue] (-1,1) circle[radius=1.6pt] node[above left] {\(i-1\)};
  \fill[polygonblue] (-2,0) circle[radius=1.6pt] node[above left=1pt] {\(-2\)};
  \fill[polygonblue] (-1,-1) circle[radius=1.6pt] node[below left] {\(-i-1\)};
  \fill[polygonrust] (4,0) circle[radius=1.8pt] node[above] {\(4\)};
  \node[polygonrust,below] at (4,-.07) {Fourth coordinate};
 \end{tikzpicture}
 \caption{For \(M=4\) and \(a=1\), the path in the equality case would
 give \(i-1,-2,-i-1\). To make their sum with the fourth coordinate
 zero, that coordinate would have to be \(4\). Its modulus exceeds
 those of the other three coordinates.}
 \label{fig:centering}
\end{figure}
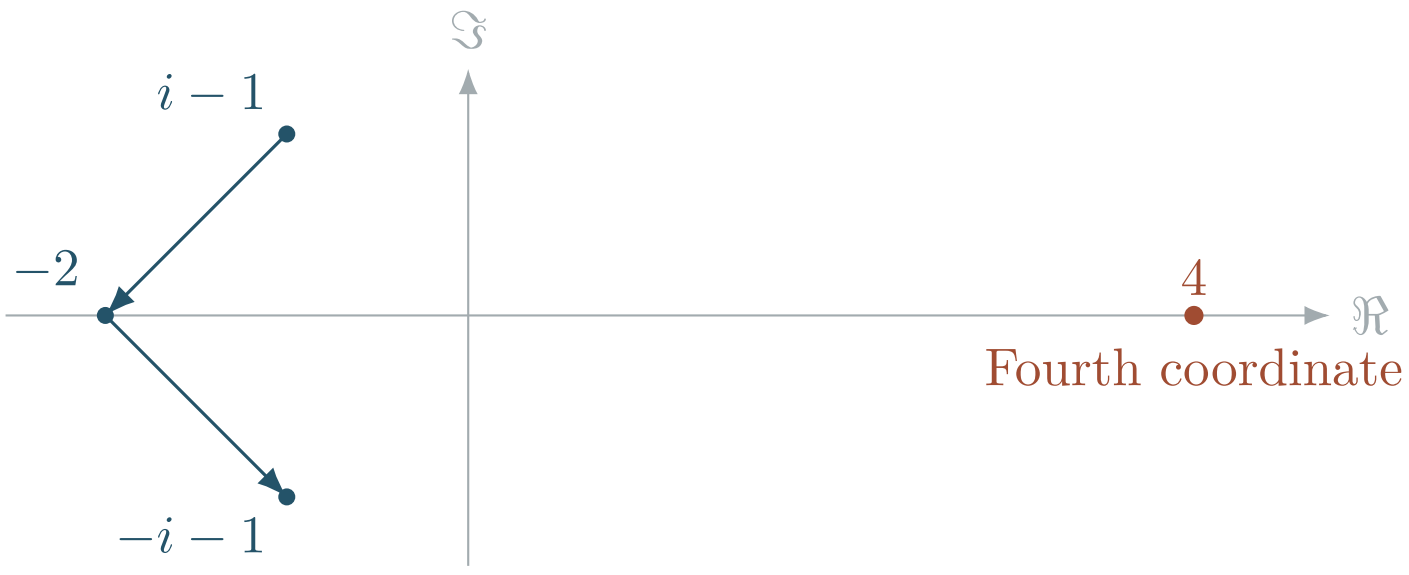

For \(M=4\) and \(a=1\), the path coordinates of the proposed equality
vector are \(i-1,-2,-i-1\). They sum to \(-4\), and the remaining
coordinate would be \(4\), as in Figure~\ref{fig:centering}. In the
example with \(\omega=-1\), the value \(4\) would have to occur together
with \(-4\). The three path coordinates already account for all the
other indices of \(S\).

\subsection{A path on all \texorpdfstring{\(M\)}{M} vertices}

When the path contains all \(M\) indices of \(S\), the added zero gives a
cycle of length \(M+1\). The spectrum of \(C_{M+1}\) contains
\(e^{2\pi i/(M+1)}\), on the open arc from \(1\) to \(\xi\) beyond
the side. On the path itself, the coordinates have sum zero by
Lemma~\ref{lem:support}. We will use this condition to estimate the
path contribution.

Let \(J_M\) be the matrix with ones immediately above its diagonal and
zeros elsewhere. Thus \(w^*J_Mw=\sum_{j=1}^{M-1}\bar w_jw_{j+1}\).
The difference between the desired upper bound and the path
contribution is \(w^*Hw\), where
\begin{equation}\label{eq:H}
 H=\cos\alpha\,I-\tfrac12(e^{-i\alpha}J_M+e^{i\alpha}J_M^*) .
\end{equation}
Thus we have to prove \(w^*Hw>0\) for nonzero \(w\) in the zero-sum
space \(V=\{w\in\mathbb C^M:\mathbf 1^*w=0\}\), where \(\mathbf 1\)
is the vector of \(M\) ones.

We start by counting the signs of the eigenvalues of \(H\) on the full
space. Let \(D\) be diagonal with \(D_{jj}=e^{ij\alpha}\). Writing
\(w=Dy\) cancels the phases in \(e^{-i\alpha}\bar w_jw_{j+1}\), so
\(D^*HD=\cos\alpha\,I-(J_M+J_M^*)/2\). The new coordinates serve
to compute the spectrum. The matrix \((J_M+J_M^*)/2\) averages two
neighbouring entries with \(y_0=y_{M+1}=0\). For
\(\beta=k\pi/(M+1)\), \(1\le k\le M\), the vector with entries
\(\sin(j\beta)\) vanishes at both ends and satisfies
\[
 \tfrac12\big(\sin((j-1)\beta)+\sin((j+1)\beta)\big)
     =\cos\beta\,\sin(j\beta).
\]
The \(M\) distinct values \(\cos\beta\) give the whole spectrum.
The eigenvalues of \(H\) are therefore
\[
 \cos\alpha-\cos\!\left(\frac{k\pi}{M+1}\right),\qquad 1\le k\le M.
\]
Since \(\pi/(M+1)<\alpha<2\pi/(M+1)\), exactly one is negative and the
others are positive. For example, at \(M=4\) the unrestricted path
contribution can reach \(\cos(\pi/5)\|w\|_2^2\), which is greater
than the required \(\cos(\pi/4)\|w\|_2^2\).

We return to \(V\) in the original coordinates. We seek a negative
line through a vector \(x\) with \(x^*Hw=0\) for \(w\in V\).
As \(H\) is Hermitian, \(x^*Hw=(Hx)^*w\), so \(Hx\) has to be
perpendicular to \(V\). The ordinary orthogonal complement of \(V\)
is spanned by \(\mathbf 1\),
so we solve \(Hx=\mathbf 1\), using the invertibility of \(H\). Now
\[
 x^*Hw=(Hx)^*w=\mathbf 1^*w .
\]
If \(\mathbf 1^*x<0\), then \(x^*Hx=\mathbf 1^*x<0\). Since
\(x\) has nonzero coordinate sum, each vector splits uniquely
as a multiple of \(x\) plus a vector in \(V\). To find the
multiple, match the coordinate sum.

If some vector of \(V\) also had a negative value, it and \(x\)
would span a two-dimensional subspace on which every nonzero
vector has a negative form value because the mixed terms vanish.
A subspace of dimension two meets the span of the \(M-1\)
eigenvectors for positive eigenvalues in a nonzero vector,
where the form is positive. Thus the form is nonnegative on \(V\).

In the basis made of \(x\) and a basis of \(V\), the matrix of the
form is block diagonal, with blocks \(\mathbf 1^*x\) and the form
on \(V\). The full matrix remains invertible under the change of
basis, so the block on \(V\) is invertible. Nonnegativity therefore
gives positive definiteness on \(V\).

It remains to find the sign of \(\mathbf 1^*x\). Extend \(x\) by
\(x_0=x_{M+1}=0\). The \(j\)th equation of \(Hx=\mathbf 1\) reads
\begin{equation}\label{eq:path-recurrence}
 \cos\alpha\,x_j
  -\tfrac12(e^{-i\alpha}x_{j+1}+e^{i\alpha}x_{j-1})=1
 \qquad(1\le j\le M).
\end{equation}
Before imposing the endpoint conditions, the same three coefficients
on the left occur in each row and sum to zero. A sequence linear in \(j\)
therefore gives a constant. For the sequence \(j\) it is \(i\sin\alpha\),
so \(j/(i\sin\alpha)\) gives the right-hand side \(1\). The sequences \(1\) and \(\xi^j\) both
give zero on the left. The difference \(\xi^j-1\) vanishes at \(j=0\)
and equals \(\xi-1\) at \(j=M+1\). A multiple of that difference
cancels the value \((M+1)/(i\sin\alpha)\) at \(M+1\), giving
\[
 x_j=\frac{1}{i\sin\alpha}
          \left(j-(M+1)\frac{\xi^j-1}{\xi-1}\right),\qquad 1\le j\le M .
\]
Summing uses \(\sum_{j=1}^M j=M(M+1)/2\) and \(\sum_{j=1}^M\xi^j=0\).
Since \((\xi+1)/(\xi-1)=-i\cot\alpha\), this gives
\begin{equation}\label{eq:negative}
 \mathbf 1^*x
   =\frac{M(M+1)}{i\sin\alpha}
     \left(\frac12+\frac1{\xi-1}\right)
   =-\frac{M(M+1)\cos\alpha}{2\sin^2\alpha}<0 .
\end{equation}

\Needspace{6\baselineskip}
In the four-coordinate case the formula gives
\[
 x=\sqrt2\,(4i,\,-5+3i,\,-5-3i,\,-4i)^\top,\qquad
 \mathbf 1^*x=-10\sqrt2 .
\]
By \eqref{eq:negative} the line through \(x\) is negative, and hence
the form is positive definite on \(V\). The path vector is nonzero
and lies in \(V\), so \eqref{eq:H} gives a strict inequality for its
contribution. Equality in \(F_\eta(p,v)\le\cos\alpha\,W_S(v)\)
requires every component to be a cycle and to attain its cycle bound.

\section{The side inequalities for nearby eigenpairs}\label{sec:limit}

For a permutation with a strict side bound at \(v\), the bound
holds near \(v\) by continuity. If a permutation attains equality,
we need more information. A plane curve can leave the half-plane
even if it is tangent to the side. For example,
\(z(t)=(1-t)\omega+t\omega\xi+t^2\eta\), with
\(\eta=\omega e^{i\alpha}\), has derivative along the side at \(0\)
and projection \(\cos\alpha+t^2\). For the equality permutations
we obtain a bound from their action on \(S\).

\begin{lemma}\label{lem:stable}
Let \(\omega\) be a root of unity of exact order \(q\ge2\), and let \(M\ge3\) be a
multiple of \(q\). Suppose \(P_rv=\omega v\) with \(v\ne0\) and
\(|\suppv v|\le M\). There is a neighbourhood of \(v\) such that
if a doubly stochastic matrix \(A\) has eigenvalue \(z\) with an
eigenvector \(u\) in that neighbourhood, then \(z\) satisfies the two side
inequalities \eqref{eq:side} of \(\Pi_M\) at \(\omega\).
\end{lemma}

\begin{proof}
Let \(S=\suppv v\), \(\alpha=\pi/M\), \(\xi=e^{2i\alpha}\) and
\(\eta=\omega e^{i\alpha}\). Suppose \(p\) attains the side bound
at \(v\). In \eqref{eq:relative} each component attains its own
bound and is a cycle, by Section~\ref{sec:paths}. This means
\(T=r^{-1}\circ p\) maps \(S\) onto itself. Since \(r\) also
preserves \(S\), the same is true of \(p\). By
Section~\ref{sec:cycle-comparison}, either the values are
constant on every cycle of \(T\), or a single cycle of length
\(M\) with values in geometric progression uses all of \(S\).

If the coordinates of \(v\) are constant around each cycle of \(T\) in
\(S\), then \(v_{p(i)}=\omega v_i\) on \(S\). A cycle of \(p\) in \(S\)
of length \(\ell\) therefore has \(q\mid\ell\), by following the
nonzero coordinate once around that cycle. The block of its matrix
has the \(\ell\)th roots of unity as eigenvalues. These include
\(\omega\), and since \(\ell\le M\), their spacing is at least
\(2\pi/M\). They all satisfy the side inequality at \(\omega\) with
normal \(\eta\).

If one cycle of \(T\) of length \(M\) uses \(S\) with
\(v_{T(i)}=\xi v_i\), then \(v_{p(i)}=\omega\xi v_i\).
Because \(q\mid M\), applying \(p\) \(M\) times multiplies
the value by \((\omega\xi)^M=1\). The values on \(S\) are
distinct, so this returns each index to itself. On \(S\)
the permutation matrix therefore has \(M\)th power the
identity and its eigenvalues are \(M\)th roots, all satisfying
the side bound.

The principal submatrix \(P_p[S,S]\) is a permutation matrix in
both cases, so its numerical range is the convex hull of its
eigenvalues and lies in the side half-plane. For \(u_S\) the
coordinates of \(u\) in \(S\), we have
\(F_\eta(p,u)=\Re(\bar\eta\,u_S^*P_p[S,S]u_S)\). The numerical
range bound gives \(F_\eta(p,u)\le\cos\alpha\,W_S(u)\) for
every \(u\).

For each \(p\) with strict inequality at \(v\), continuity gives
the bound on a neighbourhood of \(v\). Conjugating \(\omega\) and
\(v\) gives the bound for \(\omega e^{-i\alpha}\), since the
permutation matrices are real. There are finitely many \(p\)
and two normals. Intersect their neighbourhoods and shrink the
intersection so that \(W_S(u)>0\). Dividing \eqref{eq:scorebound} by that weight
gives the two bounds on \(z\).
\end{proof}

In Example~\ref{ex:square}, use the side from \(-1\) to \(-i\), so
\(\xi=i\). For \(p=r=(1\,3)(2\,4)\) fixing \(5\),
\(T\) is the identity and the cycles each consist of one index.
For \(P_p=C_4^3\oplus[1]\), \(T\) is the four-cycle of \(C_4\) on
\(S\), with coordinates \(1,i,-1,-i\) in cycle order. Varying
the convex coefficient between these two matrices traces the
side at which they give equality.

\begin{proof}[Proof of Theorem~\ref{thm:local}]
We already have \(\Pi_M\subseteq U_N\subseteq D_N\). The other
inclusion is local and we prove it using Lemma~\ref{lem:stable}.

Take \(q\ge2\) and \(M\ge3\). If every disk about \(\omega\) contained
a point of \(D_N\) outside \(\Pi_M\), take such points \(z_j\to\omega\)
and corresponding eigenpairs \(A_ju^{(j)}=z_ju^{(j)}\) with
\(\|u^{(j)}\|_2=1\). Along a subsequence \(A_j\to A\) and
\(u^{(j)}\to v\) as in Section~\ref{sec:support}.
Lemma~\ref{lem:support} gives \(P_rv=\omega v\) for some \(r\) and
\(q\mid|\suppv v|\). The size bound \(|\suppv v|\le N\) then gives
\(|\suppv v|\le M\). Lemma~\ref{lem:stable} applies with this \(M\)
and gives the two side inequalities at \(\omega\) for large \(j\).
By Section~\ref{sec:polygons}, these imply \(z_j\in\Pi_M\) for
large \(j\), a contradiction.

At \(\omega=1\) with \(N\ge3\), we have \(M=N\). The two sides of
\(\Pi_N\) at \(1\) cut off the open arcs from \(1\) to \(e^{2\pi i/N}\)
and to \(e^{-2\pi i/N}\). Every root of order at most \(N\) satisfies
both side inequalities, since every such root other than \(1\)
has angular distance at least \(2\pi/N\) from \(1\). So the
inequalities hold on \(\conv U_N\), which contains \(D_N\).

For \(N=1\) the spectrum is \(\{1\}\). In order two the diagonal
entries are both \(a\) and the off-diagonal entries both \(1-a\)
for \(0\le a\le1\). The vectors \((1,1)^\top\) and \((1,-1)^\top\)
give eigenvalues \(1\) and \(2a-1\), so \(D_2=\Pi_2\).
The remaining case \(M=2\) is \(N=3\) and \(\omega=-1\). A nonreal
eigenvalue \(z\) would give the spectrum \(1,z,\bar z\). The trace
is nonnegative and equals \(1+2\Re z\). Thus the disk
\(|z+1|<1/2\), where \(\Re z<-1/2\), contains only real
eigenvalues. The unit disk bound puts them in \(\Pi_2\).

The set \(D_N\) is compact. It lies in the closed unit disk, and
if a sequence of its points converges, choose a matrix and a unit
eigenvector for each point. Convergent subsequences of the matrices
and vectors give an eigenpair at the limit, so \(D_N\) is closed.
By Lemma~\ref{lem:support}, its points on the unit circle are
among the roots of order at most \(N\). Remove from \(D_N\) the
disks of local equality at these finitely many roots. The compact
set left over has maximum modulus \(r_N<1\) when it is nonempty.
Otherwise take \(r_N=0\). A point of \(D_N\) of modulus greater
than \(r_N\) belongs to one of the disks and hence to \(U_N\).
Together with \(U_N\subseteq D_N\) this gives the annulus.
\end{proof}

\Needspace{17\baselineskip}
\section{Approaching a fixed permutation matrix}\label{sec:fixed}

For the full region, the matrix size bounds the support of a
limiting eigenvector. If the matrix approaches a fixed permutation
matrix \(P\), the cycles of \(P\) give the bound \(L\). We can use
this bound in Lemma~\ref{lem:stable} to prove the endpoint
conjecture of Johnson, Newman and Spitkovsky \cite{jns}. For
matrices, \(\|\cdot\|_2\) denotes the operator norm induced by
the Euclidean vector norm.

\begin{corollary}\label{cor:fixed}
Let \(P=P_p\) be a permutation matrix having eigenvalue \(\omega\)
of exact order \(q\). Set
\[
 L=\sum_{\substack{C\text{ a cycle of }p\\ q\mid |C|}} |C|.
\]
There exist \(\delta,\varepsilon>0\) for which every doubly stochastic
matrix \(A\) with \(\|A-P\|_2<\delta\) satisfies
\[
 \sigma(A)\cap B(\omega,\varepsilon)\subseteq\Pi_L .
\]
In particular, let \(R\) be doubly stochastic and let \(z(t)\) be a
continuous eigenvalue branch of \(A(t)=(1-t)P+tR\), \(0\le t\le1\).
If \(z(0)=\omega\), then \(z(t)\in\Pi_L\) for all sufficiently small
\(t\ge0\).
\end{corollary}

\begin{proof}
For \(q\ge2\) and \(L\ge3\), suppose the inclusion failed for
arbitrarily small \(\delta,\varepsilon\). There would be
\(A_j\to P\), \(z_j\to\omega\) with \(z_j\in\sigma(A_j)\setminus\Pi_L\).
Unit eigenvectors have a subsequential limit \(v\) with \(Pv=\omega v\).
On a cycle of length \(\ell\) we have \(v_i=\omega^\ell v_i\), so
nonzero coordinates occur only on cycles counted by \(L\). This
gives \(|\suppv v|\le L\) and \(q\mid L\). Taking \(M=L\) in
Lemma~\ref{lem:stable} gives the two side inequalities for \(z_j\)
for large \(j\). Near \(\omega\) they give \(z_j\in\Pi_L\), which
contradicts the choice of \(z_j\).

If \(q=1\), then \(L=N\) and the assertion follows from
Theorem~\ref{thm:local}. If \(L=q=2\), the only even cycle of \(p\)
has length two, so \(-1\) is a simple eigenvalue of \(P\).
Choose a disk centred at \(-1\) whose closure is disjoint from
the other eigenvalues of \(P\). By continuity of the roots of
the characteristic polynomial, the disk continues to contain
just one eigenvalue, counted with multiplicity, for matrices
sufficiently close to \(P\). Conjugation fixes this eigenvalue
since \(A\) is real and the disk is centred at \(-1\).
It lies in \([-1,1]=\Pi_2\) by the unit disk bound.

Along the given branch \(z(t)\to\omega\). Also
\(\|A(t)-P\|_2=t\|R-P\|_2\le2t\), because \(R\) and \(P\) both
have operator norm at most one. For small \(t\), both bounds of
the neighbourhood assertion apply.
\end{proof}

\begin{figure}[!htbp]
 \centering
 \begin{tikzpicture}[x=4cm,y=4cm,>=Latex,font=\small]
  \pgfmathsetmacro{\octagonslope}{sqrt(2)-1}
  \fill[polygonfill!55] (0,0)--({.50*\octagonslope},.50)--(.64,.50)
      --(.64,-.50)--({.50*\octagonslope},-.50)--cycle;
  \fill[polygonblue!19] (0,0)--(.50,.50)--(.64,.50)
      --(.64,-.50)--(.50,-.50)--cycle;
  \draw[guidegrey!75,thin,->] (-.14,0)--(.74,0) node[right] {\(x\)};
  \draw[guidegrey!75,thin,->] (0,-.58)--(0,.59) node[above] {\(y\)};
  \draw[polygonrust,semithick,dashed] ({.50*\octagonslope},.50)--(0,0)--({.50*\octagonslope},-.50);
  \draw[polygonblue,thick] (.50,.50)--(0,0)--(.50,-.50);
  \fill[polygonblue] (0,0) circle[radius=1.6pt] node[below left] {\(z=-1\)};
  \draw[polygonrust,thin] ({.39*\octagonslope},.39)--(.25,.55)--(.65,.55);
  \node[polygonrust,right] at (.65,.55) {Octagon side};
  \draw[polygonblue,thin] (.24,-.24)--(.35,-.28)--(.65,-.28);
  \node[polygonblue,right] at (.65,-.28) {Square side};
 \end{tikzpicture}
 \caption{The side bounds from Example~\ref{ex:fixed} in the inward and transverse coordinates of
 \eqref{eq:wedge}. The square gives the narrower angle near the fixed permutation. The octagon gives the
 angle for the full region \(D_8\).}
 \label{fig:fixed}
\end{figure}
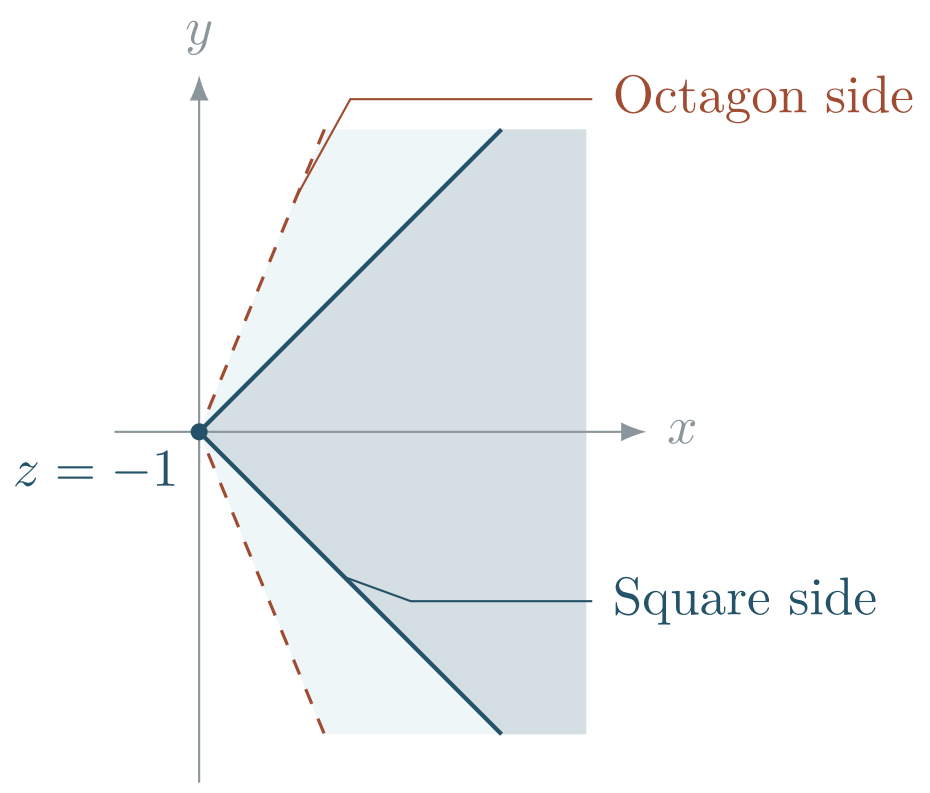

\begin{example}[Two cycles of length two in order eight]\label{ex:fixed}
Take
\[
 P=C_4^2\oplus C_3\oplus[1],\qquad \omega=-1.
\]
Here \(M=8\) and \(L=4\), since only the two cycles of length two
among \(2,2,3,1\) contribute to \(L\). Corollary~\ref{cor:fixed}
therefore gives the sides of the square at \(-1\) for matrices
near \(P\). In the coordinates of \eqref{eq:wedge} they give
\(x\ge|y|\). Near \(-1\), the full region has the sides of
the octagon, which give \(x\ge|y|\tan(\pi/8)\).

We can attain both sides of the square using the construction
in Example~\ref{ex:square}. Put \(R_+=C_4\oplus C_3\oplus[1]\)
and \(R_-=C_4^3\oplus C_3\oplus[1]\). For \(A_\pm(t)=(1-t)P+tR_\pm\),
the vector \(v=(1,i,-1,-i,0,0,0,0)^\top\) gives the
eigenvalues \(-1+t(1\pm i)\). The blocks \(C_3\) and \([1]\)
act on zeros, so they leave the calculation unchanged.

\Needspace{4\baselineskip}
For \(0<t\le1\), take convex combinations of \(P,A_+(t)\) and
\(A_-(t)\). The matrices stay within distance \(2t\) of \(P\),
and their eigenvalues on \(v\) fill the triangle spanned by
\(-1,-1+t(1+i),-1+t(1-i)\). This triangle contains the part of
the square sufficiently close to \(-1\). As \(t\) can be
arbitrarily small, the local square bound is attained at \(P\).
Figure~\ref{fig:fixed} compares the two angles.
\end{example}

\smallskip\noindent
The local equality and the annulus assertion of Theorem~\ref{thm:local}
appeared in Theorem~1.2 of
\href{https://arxiv.org/abs/2611.30905v1}{arXiv:2611.30905v1}.
This paper develops that proof and includes the fixed-permutation bound
in Corollary~\ref{cor:fixed}.

\smallskip\noindent\textbf{Use of AI tools.}
The proof was developed through conversations with ChatGPT~6 Max in Codex.
Codex also helped draft and revise the exposition, prepare the TikZ figures,
and run a numerical check of recurrence~\eqref{eq:path-recurrence}. The authors take full
responsibility for the mathematical claims, proofs, computations, figures,
and final text.

\end{document}